\documentclass[a4paper,10pt]{ifacconf}
\usepackage{graphicx}
\usepackage{amssymb}
\usepackage{amsmath}
\usepackage{epstopdf}
\usepackage{natbib}   
\usepackage{csquotes}
\usepackage{algorithm}
\usepackage[noend]{algpseudocode}

\newtheorem{remark}{Remark}
\newtheorem{definition}{Definition}

\begin{document}

\begin{frontmatter}

\title{Information Structure Design in\\ Team Decision Problems}

\thanks[footnoteinfo]{The work of T. Summers is supported by the National Science Foundation under Grant CNS-1566127. M. Kamgarpour is supported by the European Union ERC Starting Grant, CONENE. e-mail: tyler.summers@utdallas.edu, mkamgar@control.ee.ethz.ch.}

\author{Tyler Summers$^\dag$ \quad } 
\author{Changyuan Li$^\dag$ \quad }
\author{Maryam Kamgarpour$^\ddag$} 

\address{$^\dag$University of Texas at Dallas \quad $^\ddag$ETH Z\"urich}

\begin{abstract}                
We consider a problem of information structure design in team decision problems and team games. We propose simple, scalable greedy algorithms for adding a set of extra information links to optimize team performance and resilience to non-cooperative and adversarial agents. We show via a simple counterexample that the set function mapping additional information links to team performance is in general not supermodular. Although this implies that the greedy algorithm is not accompanied by worst-case performance guarantees, we illustrate through numerical experiments that it can produce effective and often optimal or near optimal information structure modifications. 
\end{abstract}


\end{frontmatter}

\section{Introduction}
Future critical infrastructures, including electric power, transportation, water, etc., are emerging as cyber-physical networks that will feature cooperative autonomous decision making agents equipped with embedded sensing, computation, communication, and actuation capabilities. A key challenge in the analysis and design of these networks is decentralization of information: each decision making agent must act based on partial information measured or received locally to optimize network operation. Moreover, with growing concerns over cyber-physical security \cite{cardenas2008secure,zhu2011taxonomy,zhu2011hierarchical,pasqualetti2013attack,sandberg2015cyberphysical}, each agent must not only coordinate its actions with team members, but must also counter against teams of malicious agents to mitigate attack impacts and provide resiliency. The information structure -- who knows what and when -- is a basic component in formal analyses of these issues and plays a crucial role in determining optimal strategies and computational tractability \cite{radner1962team,witsenhausen1971information,ho1972team,basar1978decentralized,ho1980team,rotkowitz2006characterization,nayyar2013decentralized,yuksel2013stochastic,bacsar2014stochastic,lessard2015optimal}.

While the importance of information structure is widely recognized in team decision theory, decentralized control, and game theory, the vast majority of the literature focuses on designing decision and control strategies for a \emph{given} information structure. The \emph{design} of information structures -- who \emph{should} know what and when -- has been recognized as an important problem since the earliest work on team decision theory \cite{marschak1955elements,radner1962team}, but has received very little formal attention. Radner notes the emphasis on analysis of strategies for given information structures in a seminal paper on team decision theory:
\begin{center}
\begin{displayquote}
\emph{An important organizational problem is the determination of what statistical information shall be made available to the various decision makers in the organization. Implicit in the solution of such a problem is the determination of the best use that can be made of any given structure of information, i.e., the best decision functions. The results to be presented here are concerned with this latter problem.}
\end{displayquote}
\end{center}
This emphasis on analysis of given fixed information structure has followed in much of the related work, and has presented rich challenges for many decades.

We believe it is important to shift some focus to information structure design, where one \emph{jointly} optimizes the information structure together with decision strategies, especially in the context of emerging cyber-physical networks. We consider problems of information structure design in team decision problems and team games. We focus here on static problems since no work to our knowledge has been done even in this setting.
We also focus on linear quadratic problems since they are analytically tractable, admitting closed-form equilibrium solutions that provide insight into essential properties. In the non-cooperative game setting, we focus on a specific class of games involving two teams with decentralized information structure \cite{colombino2015quadratic}, maintaining a sharp distinction between cooperative and adversarial features. 

The main contributions are as follows. We formulate several information structure design problems as (finite, combinatorial) set function optimization problems. These can be solved in principle by brute force enumeration, but this approach is not feasible even for moderately sized networks, and is certainly ineffective for the large networks of critical infrastructure that motivate this work. We therefore propose simple greedy algorithms that provide an effective and scalable heuristic. We show via a simple counterexample that the set function mapping additional information links to team performance is in general not supermodular. Although this implies that the greedy algorithm is not accompanied by theoretical worst-case performance guarantees, we illustrate its effectiveness and scalability through numerical experiments, showing that it often produces optimal or near optimal information structure modifications. 

Our focus here is on a general mathematical framework, but many emerging applications in cyber-physical networks feature distributed estimation and control problems that can be formulated as team decision problems and games. For example, information structure design in electrical power networks includes optimal sensor placement (e.g., phasor measurement units and other advanced metering) and optimal communication design for wide area monitoring and control. Furthermore, many large interconnected power grids are operated by a set of independent transmission system operators  with objective functions that are not necessarily aligned, and are susceptible to influence by distributed attacking teams with adversarial objectives. Similar distributed estimation and control problems can be formulated in other critical infrastructure, such as transportation networks. 



The rest of the paper is laid out as follows. In Section 2 we formulate information structure design problems for single team and multiple team decision problems. In Section 3 we present a greedy algorithm for information structure design. Section 4 presents numerical experiments. Section 5 gives concluding remarks. 

%
%
%
%
%
%

\section{Problem formulation}
We formulate two separate information structure design problems. The first is a (cooperative) single team problem, and the second is a two team problem, which has \emph{both} cooperative and non-cooperative/adversarial features.

\subsection{Team Decision Problems}

\paragraph*{Fixed information structure.} A team decision problem involves coordinating the decisions of a team of $N$ decision making agents in a stochastic environment.  The state of the environment is assumed to be a normal random vector $x \in \mathbf{R}^n$ with mean $\bar x = \mathbf{E} x$ and covariance matrix $X = \mathbf{E} x x^T \succ 0$, i.e., $x \sim \mathcal{N}(\bar{x}, X)$. It is assumed that every agent knows the environment state statistics $\bar x$ and $X$. In addition, each agent of the team receives its own noisy local information about the environment state $x$, which we assume to be linear:
\begin{equation}
z_i = H_i x + w_i, \quad i = 1,.., N
\end{equation}
where $H_i \in \mathbf{R}^{p_i \times n}$ and $w_i \sim \mathcal{N}(0,R_i)$ with $R_i \succeq 0$. Each row of $H_i$ can represent information obtained from a sensor or a communication link with another device or agent in the team. For example, in a power network each agent may be a local network monitoring station, and $H_i$ could include local measurements from phasor measurement units or communicated information from other parts of the network. We define the information structure as a collection of the parameters specifying the information for each agent $$S_0 = \{ (H_1, R_1), (H_2, R_2), ..., (H_N, R_N) \}.$$

Each agent must select a decision function $\gamma_i : \mathbf{R}^{p_i} \rightarrow \mathbf{R}^{m_i}$ from a set of Borel measurable functions that specifies its decision $u_i = \gamma_i(z_i)$ based on realizations of the random variable $z_i$. We define the team decision function $\gamma = (\gamma_1,...,\gamma_N)$ and the associated team decision vector $u = [u_1^T,...,u_N^T]^T \in \mathbf{R}^{\Sigma_i m_i}$, which may represent a distributed parameter estimate or control action. The quadratic team cost function is
\begin{equation}
\bar J(u) = u^T Q x + \frac{1}{2} u^T P u,
\end{equation}
where $$Q = \left[\begin{array}{c}Q_{1} \\Q_{2} \\\vdots \\Q_{N}\end{array}\right], \quad P = \left[\begin{array}{cccc}P_{11} & P_{12} & \cdots & P_{1N} \\P_{12}^T & P_{22} & \cdots & P_{2N} \\\vdots & \vdots & \ddots & \vdots \\P_{1N}^T & P_{2N}^T & \cdots & P_{NN}\end{array}\right] $$
with $Q$ and $P$ partitioned according to agent decision dimensions, i.e., $Q_i \in \mathbf{R}^{m_i \times n}$, $P_{ii} \in \mathbf{R}^{m_i \times m_i}$, and $P_{ij} \in \mathbf{R}^{m_i \times m_j}$. The cost function is assumed to be strictly convex in $u$, i.e., $P \succ 0$. For any given team decision function $\gamma$ we define the expected cost
\begin{equation}
J(\gamma, S_0) = \mathbf{E}_x (u^T Q x + u^T P u), \quad u_i = \gamma_i(z_i(x))
\end{equation}
The optimal value of the objective function under the optimal team decision function is denoted by
\begin{equation}
J^*(S_0) = \min_\gamma J( \gamma ), \quad \gamma^* = \arg \min_\gamma J(\gamma)
\end{equation}

Under the stated assumptions, the optimal decision functions $\gamma_i^*$ are affine and can be computed by solving a set of linear equations derived from stationarity conditions; see \cite{radner1962team}, or for a more general multi-objective game formulation \cite{basar1978decentralized}. In particular, the optimal solution consists of each agent forming the conditional state estimate
\begin{equation}
\begin{aligned}
\hat x_i &= \mathbf{E}[x \mid z_i]  \\
             &= \bar x + X H_i^T (H_iXH_i^T + R_i)^{-1} (z_i - H_i \bar x)
\end{aligned}
\end{equation}
and using the affine decision function
\begin{equation}
u_i = \gamma_i(z_i) = A_i \bar x + B_i (\hat x_i - \bar x),
\end{equation}
where $A_i$ and $B_i$ are the unique solutions to the linear equations
\begin{equation} \label{teamsol}
\begin{aligned}
P_{ii} A_i +  \sum_{j \neq i} P_{ij} A_j &= -Q_i  \\
P_{ii} B_i  + \sum_{j \neq i} P_{ij} B_j X H_j^T (H_jXH_j^T + R_j)^{-1} H_j &= -Q_i.
\end{aligned}
\end{equation}

\paragraph*{Information Structure Design.}
We now suppose that for each agent there is a finite set of possible measurements or communicated information about the environmental state that could be added to its information; we let $q_i$ denote the number of possible additional measurements or communication links that could be added to agent $i$. We collect the parameters defining these possibilities for the whole team into the finite set
\begin{equation} \label{Vset}
\begin{aligned}
V = \{ &(h_{11},r_{11}), (h_{12},r_{12}), ..., (h_{1q_1},r_{1q_1}), \\
          &(h_{21},r_{21}), (h_{22},r_{22}), ..., (h_{2q_2},r_{2q_2}) ..., \\
          &(h_{N1},r_{N1}), (h_{N2},r_{N2}), ..., (h_{Nq_N},r_{Nq_N}) \}
\end{aligned}
\end{equation}
where $h_{ij} \in \mathbf{R}^n$ represents the $j$th possible additional measurement or communicated information about the environmental state that could be added to the information of agent $i$, and $r_{ij} \geq 0$ represents the associated variance. We assume that each additional observation has an associated measurement noise that is independent of other measurement noise variables.\footnote{It is straightforward to allow noise of additional observation to be statistically dependent on other noise variables, but we assume independence to simplify notation.} In a power network, $V$ may represent, e.g., a set of additional phasor measurement units or wide area communication links that could augment the information set of each agent.


For any subset $S \subseteq V$, we associate a modified information structure by including the selected information links in the appropriate agents' information model. 
For example, the information structure modification $$S = \{ (h_{13}, r_{13}), (h_{32},r_{32}), (h_{43},r_{43}), (h_{45},r_{45}) \} \subset V$$ means that we add the third possible additional link to agent 1, the second possible additional link to agent 3, and the third and fifth possible additional links to agent 4, so that
\begin{equation}
\begin{aligned}
 z_1 &= \left[\begin{array}{c}H_1 \\h_{13}^T\end{array}\right]x + \left[\begin{array}{c}w_1 \\w_{13}\end{array}\right],  \quad  z_3 = \left[\begin{array}{c}H_3 \\h_{32}^T\end{array}\right]x + \left[\begin{array}{c}w_3 \\w_{32}\end{array}\right], \\ 
 z_4 &= \left[\begin{array}{c}H_4 \\h_{43}^T \\ h_{45}^T\end{array}\right]x + \left[\begin{array}{c}w_4 \\w_{43} \\ w_{45}\end{array}\right],
\end{aligned}
\end{equation}
where $w_{13} \sim \mathcal{N}(0, r_{13})$, $w_{32} \sim \mathcal{N}(0, r_{32})$, $w_{43} \sim \mathcal{N}(0, r_{43})$, and $w_{45} \sim \mathcal{N}(0, r_{45})$ are independent of all other measurement noise variables.

Let $J^*(S)$ denote the optimal value of the team cost function associated with the information structure modification $S$. Our first problem of interest is to select an information structure modification of size $k$ to minimize the optimal value of the team decision problem using the associated optimal decision functions for the modified information structure\footnote{Our algorithms can be easily adapted to a setting where each information structure modification has its own fixed cost, and we search for an optimal modification that satisfies a total budget constraint.}. We can pose this as a cardinality constrained set function optimization problem 
\begin{equation} \label{teamsetfunc}
\min_{S \subset V, \ |S| = k} J^*(S). 
\end{equation}

\subsection{Two Team Games}
We now formulate an analogous problem for a two team stochastic game. In this setting, there are two teams, which we call blue and red, each of which consists of a set of decision making agents interacting in a stochastic environment. We assume again that the environment state is a normal random vector $x \in \mathbf{R}^n$ with mean $\bar x = \mathbf{E} x$ and covariance matrix $X = \mathbf{E} x x^T \succ 0$ and that every agent knows the environment state statistics $\bar x$ and $X$. The blue team has $N$ decision making agents, and the red team has $M$ decision making agents. The blue team may represent agents associated with a network operator, while the red team may represent a set of non-cooperative agents or malicious attackers. The difference here is that each team has its own objective function, introducing a non-cooperative or adversarial element to the problem in addition to the cooperation required amongst team members.
\begin{figure}
\begin{centering}
\resizebox{0.9\linewidth}{!}{\includegraphics{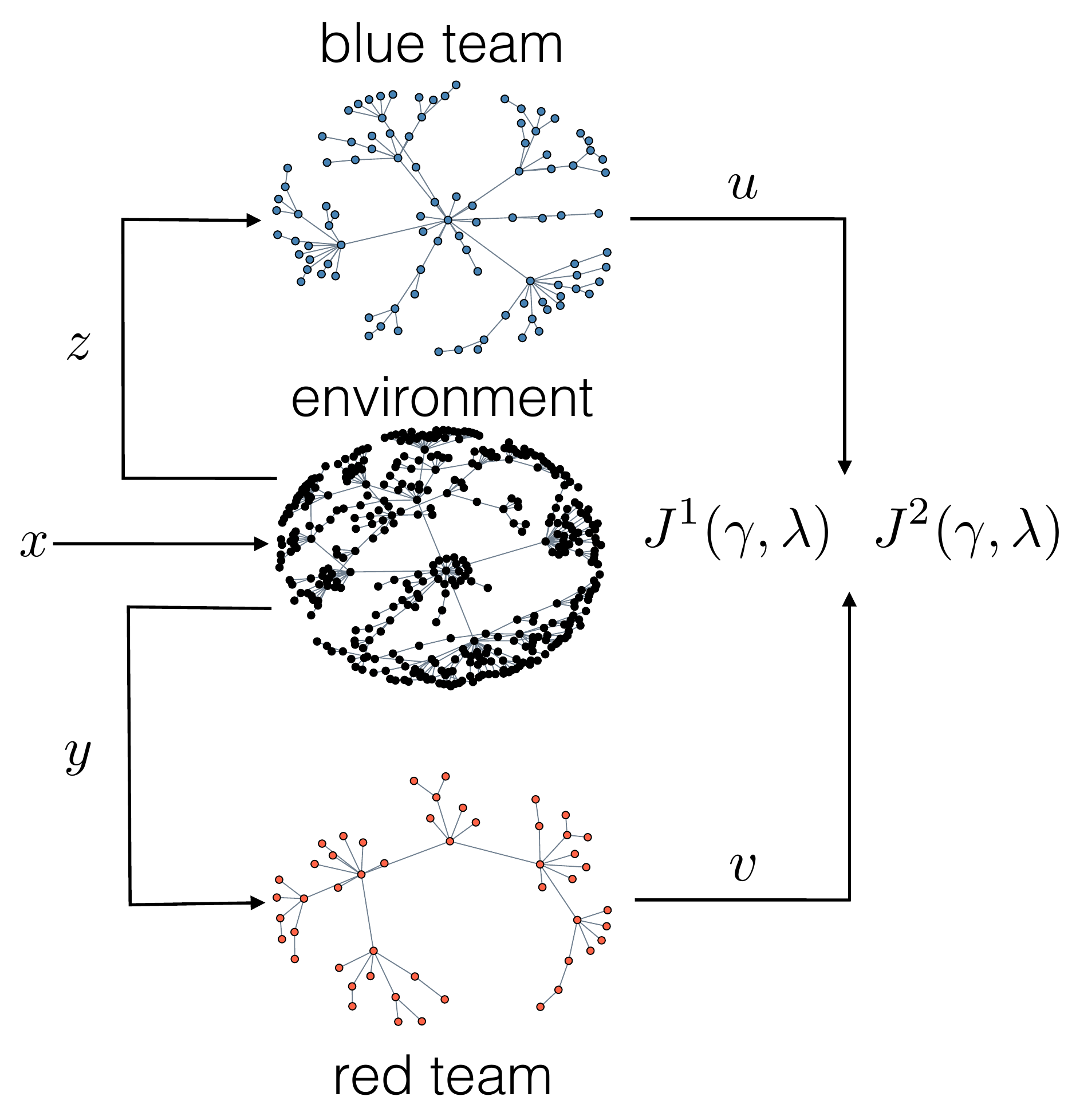}}
\caption{Illustration of the two team game setting. The blue team receives measurements $z$ and chooses a distributed decision function $u = \gamma(z)$ to optimize $J^1(\gamma,\lambda)$, and the red team receives measurements $y$ and chooses a distributed decision function $v = \lambda(y)$ to optimize $J^2(\gamma,\lambda)$. }
\end{centering}
\end{figure}

\paragraph*{Fixed information structure.}  The blue team receives information 
\begin{equation}
z_i = H_i x + w_i, \quad i = 1,.., N
\end{equation}
where $H_i \in \mathbf{R}^{p_i \times n}$ and $w_i \sim \mathcal{N}(0,R_i)$, and the red team receives information
\begin{equation}
y_j = G_j x + t_j, \quad j = 1,.., M
\end{equation}
where $G_j \in \mathbf{R}^{l_j \times n}$ and $t_j \sim \mathcal{N}(0,T_j)$. The information structure for the blue team is
$$S_0 = \{ (H_1, R_1), (H_2, R_2), ..., (H_N, R_N) \},$$
and the information structure for the red team is
$$T_0 = \{ (G_1, T_1), (G_2, T_2), ..., (G_N, T_N) \}.$$

Each agent on the blue team must select a decision function $\gamma_i : \mathbf{R}^{p_i} \rightarrow \mathbf{R}^{m_i}$ from a set of measurable functions that specifies its decision $u_i = \gamma_i(z_i)$, and each agent on the red team must select a decision function $\lambda_j : \mathbf{R}^{l_j} \rightarrow \mathbf{R}^{k_j}$ from a set of measurable functions that specifies its decision $v_j = \lambda_j(y_j)$.

In a non-cooperative two team game, each team has a separate objective function that is neither directly aligned nor misaligned with that of the opposing team. We define the team decision functions $\gamma = (\gamma_1,...,\gamma_N)$ and $\lambda = (\lambda_1,...,\lambda_M)$ and the associated team decision vectors $u = [u_1^T,...,u_N^T]^T \in \mathbf{R}^{\Sigma_i m_i}$ and $v = [v_1^T,...,v_N^T]^T \in \mathbf{R}^{\Sigma_j k_j}$. 
The blue team cost function is
\begin{equation}
\bar J^{1}(u,v) = u^T Q^{1} x + \frac{1}{2}(u^T P^{1} u + v^T R^{1} u) ,
\end{equation}
and the red team seeks to optimize a cost function
\begin{equation}
\bar J^{2}(u,v) = v^T Q^{2} x + \frac{1}{2}(v^T P^{2} v + v^T R^{2} u) ,
\end{equation}
It is assumed that $ P^{i} \succ 0, \ i=1,2$, so that $\bar J^{1}(u,v)$ is strictly convex in $u$ and $\bar J^{2}(u,v)$ is strictly convex in $v$. 

For any given team decision functions $\gamma$ and $\lambda$ we define the expected costs
\begin{equation}
\begin{aligned}
J^{1}(S_0,T_0,\gamma, \lambda) = \mathbf{E}_x (u^T Q^{1} x + u^T P^{1} u + 2v^T R^{1} u),
\end{aligned}
\end{equation}
\begin{equation}
\begin{aligned}
J^{2}(S_0,T_0,\gamma, \lambda) = \mathbf{E}_x (v^T Q^{2} x + v^T P^2 v + 2v^T R^2 u), \\ \quad \text{with} \quad u_i = \gamma_i(z_i(x)), v_j = \lambda_i(y_j(x))
\end{aligned}
\end{equation}
A pair of team decision strategies $(\gamma^*, \lambda^*)$ are called Nash equilibrium strategies if
\begin{equation}
\begin{aligned}
\gamma^* \in \arg \min_\gamma J^{1}(S_0,T_0,\gamma, \lambda^*) \\
\lambda^* \in \arg \min_\lambda J^{2}(S_0,T_0,\gamma^*, \lambda),
\end{aligned}
\end{equation}
and the corresponding Nash equilibrium values are denoted by 
\begin{equation}
\begin{aligned}
J^{1*}(S_0,T_0) =  J^{1}(S_0,T_0,\gamma^*, \lambda^*) \\
J^{2*}(S_0,T_0) =  J^{2}(S_0,T_0,\gamma^*, \lambda^*).
\end{aligned}
\end{equation}

Under the stated assumptions, the Nash equilibrium decision strategies $\gamma^*_i$ and $\lambda^*_j$ are unique and affine, and can be computed by solving a set of linear equations derived from stationarity conditions; see \cite{basar1978decentralized}. In particular, the Nash equilibrium solution also consists of each agent on each team forming the conditional state estimates
\begin{equation}
\begin{aligned}
\hat x^1_i &= \mathbf{E}[x \mid z_i]  \\
             &= \bar x + X H_i^T (H_iXH_i^T + R_i)^{-1} (z_i - H_i \bar x)
\end{aligned}
\end{equation}
\begin{equation}
\begin{aligned}
\hat x^2_j &= \mathbf{E}[x \mid y_j]  \\
             &= \bar x + X G_j^T (G_jXG_j^T + T_j)^{-1} (y_j - G_j \bar x)
\end{aligned}
\end{equation}
and using the affine decision functions
\begin{equation}
u_i = \gamma_i(z_i) = A_i \bar x + B_i (\hat x^1_i - \bar x),
\end{equation}
\begin{equation}
v_j = \lambda_j(y_j) = C_j \bar x + D_j (\hat x^2_j - \bar x),
\end{equation}
where $A_i$ and $B_i$ are the unique solutions to the linear equations
\begin{equation} \label{twoteamsol1}
\begin{aligned}
P^1_{ii} A_i +  \sum_{j \neq i}^N P^1_{ij} A_j &= -Q^1_i  \\
P^1_{ii} B_i  + \sum_{j = 1, j \neq i}^N P^1_{ij} B_j X H_j^T (H_jXH_j^T + R_j)^{-1} H_j \\
+ \sum_{j = 1, j \neq i}^M R^1_{ij} B_j X H_j^T (H_jXH_j^T + R_j)^{-1} H_j &= -Q^1_i,
\end{aligned}
\end{equation}
and $C_i$ and $D_i$ are the unique solutions to the linear equations
\begin{equation} \label{twoteamsol2}
\begin{aligned}
P^2_{ii} C_i +  \sum_{j \neq i}^N P^2_{ij} C_j &= -Q^2_i  \\
P^2_{ii} D_i  + \sum_{j=1,j \neq i}^N P^2_{ij} D_j X G_j^T (G_jXG_j^T + T_j)^{-1} G_j  \\ 
+ \sum_{j=1, j \neq i}^M R^2_{ij} D_j X G_j^T (G_jXG_j^T + T_j)^{-1} G_j  &= -Q^2_i,
\end{aligned}
\end{equation}
with $P^1$, $Q^1$, $R^1$, $P^2$, $Q^2$, $R^2$ partitioned according to the dimensions of $u_i$ and $v_j$.

\paragraph*{Information Structure Design.}
We now pose an information structure design problem for the blue team, with the information structure of the red team held fixed; an analogous problem can be posed for the red team. As above, we form the finite set $V$ in \eqref{Vset} consisting of all possible measurements or communicated information about the environmental state that could be added to the information structure of the blue team. For any subset $S \subseteq V$, we associate a modified information structure by including the selected information links in the appropriate agents' information model. 

Let $J^{1*}(S)$ denote the Nash equilibrium value of the blue team associated with the information structure modification $S$. Our second problem of interest is to select an information structure modification of size $k$ to minimize the Nash equilibrium value of the blue team under the associated Nash equilibrium strategies. Again, we can pose this as a cardinality constrained set function optimization problem
\begin{equation} \label{twoteamsetfunc}
\min_{S \subset V, \ |S| = k} J^{1*}(S). 
\end{equation}

\remark{In adversarial settings, resilient information structure design problems can be formulated for zero-sum games as a special case of the above by setting $J^1 = -J^2$, i.e., the blue team seeks to minimize $J^1$ while the red team seeks to maximize it.} 

\remark{One can also formulate variations where the blue or red team is allowed to modify the information structure of the other team (perhaps by adding links when the objectives are relatively aligned, or to sabotage by removing links or increasing noise when the objectives are relatively unaligned), or a meta-game where both teams are simultaneously allowed to modify their information structures.


\section{Information Structure Design and Lack of Supermodularity}
In this section we propose a simple greedy algorithm for the set function optimization problems defined above to formalize information structure design in team decision and game problems. We show that the set functions are not in general supermodular. This implies that the information structure modifications produced by the greedy algorithm are not in general guaranteed to come along with worst-case theoretical suboptimality gurantees. Nevertheless, the greedy algorithm can scale to far larger networks than exhaustive search, and we will demonstrate empirically that it often produces near optimal designs. 

\subsection{Set functions and submodularity}
The information structure problems described above are formulated as cardinality constrained set function optimization problems. These problems are combinatorial and finite, and so can be solved simply by brute force enumeration and exhaustive search. However, this approach quickly becomes intractable even for moderately sized problems. The motivating context of large cyber-physical networks requires a different approach.

%

Greedy algorithms are a simple alternative to exhaustive search. When a set function minimization problem has a certain property called \emph{supermodularity}, a greedy algorithm achieves results that are provable within a constant factor of the optimal value. Supermodularity (and the closely related submodularity)  plays a similar role in combinatorial optimization as convexity and concavity play in continuous optimization \cite{lovasz1983submodular,krause2012submodular}. 

\begin{definition} \label{submoddef}
A set function $f: 2^V \rightarrow \mathbf{R}$ is called \emph{supermodular} if for all subsets $A \subseteq B \subseteq V$ and all elements $s \notin B$, it holds that
\begin{equation} \label{submod1}
f(A \cup \{s\}) - f(A) \leq f(B \cup \{s\}) - f(B),
\end{equation}
or equivalently, if for all subsets $A,B \subseteq V$, it holds that
\begin{equation} \label{submod2}
f(A) + f(B) \leq f(A\cup B) + f(A\cap B).
\end{equation}
A set function is called \emph{submodular} if the reversed inequalities in (\ref{submod1}) and (\ref{submod2}) hold and is called \emph{modular} if  (\ref{submod1}) and (\ref{submod2}) hold with equality.
\end{definition}

Intuitively, supermodularity is a diminishing returns property where adding an element to a smaller set gives a larger benefit than adding it to a larger set. 
Minimization of supermodular functions (equivalently, maximization of submodular functions) is NP-hard, but a simple greedy heuristic can be used to obtain a solution that is provably close to the optimal solution \cite{greedybound}. The greedy algorithm for set function minimization is shown in Algorithm \ref{algorithm:greedy}. Several problems in systems and control that feature greedy algorithms and sub- or supermodularity have been recently explored \cite{bushnell2014supermodular,bushnell2014minimizing,summers2014submodularity,summers2014optimal,shames2014rigid,tzoumas2015sensor}. However, other important set function optimization problems in systems and control fail to be sub- or supermodular \cite{summers2016actuator}.


\begin{algorithm}[!h] 
\caption{A greedy algorithm for set function optimization.}
\label{algorithm:greedy}
\begin{algorithmic}
\State $S \leftarrow \emptyset$
\While {$|S| \leq k$}
\State $e^\star = \mathop{\mbox{argmin}}\limits_{e\in V\setminus S } \quad f(S\cup\{e\})$
\State $S \leftarrow S \cup \{e^\star\}$
\EndWhile
\State $S^\star \leftarrow S$
\end{algorithmic}
\end{algorithm}


\subsection{A greedy algorithm and lack of supermodularity}
The simple greedy algorithm described in Algorithm \ref{algorithm:greedy} can be directly applied to the information structure design problems that we formulated as cardinality constrained set function optimization problems in \eqref{teamsetfunc} and \eqref{twoteamsetfunc}. At each iteration, one simply adds the information link that reduces the optimal cost the most by evaluating the optimal cost associated with each possible additional link. The algorithm terminates after $k$ links have been added.

For the single team information structure design problem, each iteration requires a set of $2n \sum_i m_i$  linear equations to be solved to compute the $2n \sum_i m_i$ optimal strategy coefficients $A_i$ and $B_i$ in \eqref{teamsol}, so the total computational complexity is order $k|V|(n \sum_i m_i)^3$. Within each iteration, the function evaluations for computing the cost of each possible additional link are trivially parallelizable, so distributed computing platforms could be used to scale computations to large networks. Further, it may be possible to exploit the sparsity often found in many cyber-physical networks that motivate these problems.

Unfortunately, it turns out that the set functions defined in \eqref{teamsetfunc} and \eqref{twoteamsetfunc} that map information structure modifications to associated optimal team cost values or Nash equilibrium values are not in general supermodular. Consider a single team (cooperative) problem with 2 players, each of whose information could be modified by a single additional link. Suppose\\
	\begin{displaymath}
		Q=\left[ 
		\begin{array}{cc}
			Q_{1}\\
			Q_{2}
		\end{array} 
		\right]=\left[ 
		\begin{array}{cc}
			1 & 1\\
			1 & 1
		\end{array} 
		\right], \quad 
		P=\left[ 
		\begin{array}{cc}
			1 & -0.5\\
			-0.5 & 1
		\end{array} 
		\right]
	\end{displaymath}
	\begin{displaymath}
		H_{1}=H_{2}=\left[
		\begin{array}{cc}
			1 & 0
		\end{array}
		\right], \quad
		R_{1}=R_{2}=0, \quad 
		h_{11}=h_{21}=\left[
		\begin{array}{cc}
			0 & 1
		\end{array}
		\right]
	\end{displaymath}
	\begin{displaymath}
	        r_{11}=r_{21}=0, \quad
		\bar{x}=\left[
		\begin{array}{cc}
			0  \\
			0
		\end{array}
		\right], \quad
		X=\left[
		\begin{array}{cc}
			1 & 0\\
			0 & 1
		\end{array}
		\right]	
	\end{displaymath}
%

Let $V =  \{(h_{11},r_{11}), ( h_{21},r_{21} ) \}$, which has four subsets: $A=\{(h_{11},r_{11})\}$, $B=\{ ( h_{21},r_{21} ) \}$, $A\cap B = \emptyset$, and $A \cup B = V$. Evaluating the cost of all of these information structure modifications, we have:
\begin{equation} \nonumber
\begin{aligned}
J^*(\emptyset) = -2, \ J^*(A) = -2.5, \ J^*(B) = -2.5, \  J^*(V) = -4
\end{aligned}
\end{equation}
so that
$$J^*(A)+J^*(B)-J^*(A \cup B)-J^*(A \cap B)= 1 > 0$$
which violates the supermodularity inequality in Definition \ref{submoddef}. Effectively, the cost benefit provided by each additional link individually is less than the benefit of adding both of them together, so there is no diminishing returns property. It is also easy to construct examples where the submodularity inequality does not hold, so that the set function is in general neither sub- nor supermodular.  Since the single team is a special case of the two team problem, the set function for the Nash equilibrium value is neither sub- nor supermodular.

This implies that the greedy algorithm does not in general produce information structure modifications that are within a constant factor of the globally optimal information structure modifications of a given cardinality. However, the greedy algorithm can be an effective and scalable heuristic, which we demonstrate empirically next.


\section{Numerical experiments}
To illustrate the effectiveness of our proposed greedy algorithms for information structure design, we considered problems with randomly generated data that were small enough to solve globally by exhaustive search. The data was generated in the following way.
We consider a single team (cooperative)  problem with 10 states and 4 players. Each player has 3 decision variables and 2 measurements of the state, and the set of information structure modifications consists of 2 possible additional measurements for each player. The goal is to find the 4 best new measurements (out of the 8 possible) to minimize the team cost function. We let $P=\tilde{P}^{T}\tilde{P}, X=\tilde{X}^{T}\tilde{X}$ to ensure that $P$ and $X$ are symmetric and positive definite, while each element of $Q, \tilde{P},H_{i},R_{i},h_{ij},r_{ij},\tilde{X}$ are generated independently from a standard normal distribution $\mathcal{N}(0,1)$.

We compare the information structure obtained by the greedy algorithm with the globally optimal information structure found by exhaustive search. For this problem size, the greedy algorithm is about 60 times faster than exhaustive search. We observe that the greedy algorithm often finds a structure with the same value as or very near the globally optimal value. In several hundred problem instances, the greedy algorithm achieves the globally optimal value around 80\% of the time, while in the worst instance is only 25\% worse than the globally optimal value. Although there are no guarantees,  our experiment shows that the greedy algorithm can produce very good results. Moreover, it scales to problem sizes far larger than what can be handled by exhaustive search, making it much more suitable for scaling to problems involving distributed estimation and control in large cyber-physical networks. 

\section{Conclusions and Outlook}
\vspace{-0.2cm}
We have formulated information structure design problem for team decision problems and team games, in which the objective is to \emph{jointly} design information structure modifications together with optimal strategies. We posed these as set function optimization problems and proposed a greedy algorithm as a heuristic for designing good information structures. We showed via a simple counterexample that the associated set functions are in general not supermodular, so that the greedy algorithms do not in general come with worst-case performance guarantees. However, we observed empirically that the greedy algorithm often produces effective information structure modifications.

Our immediate future work will consider team decision problems and games with dynamics, focusing on tractable information structures in that setting, such as partially nested and quadratically invariant. We will also explore alternative convex relaxation approaches and other techniques for scaling the computations to large networks. Finally, we plan to apply the results to application areas, including power systems and transportation networks.

\vspace{-0.1cm}

\bibliography{refs}  

\end{document}